\documentclass[12pt,notitlepage,twoside]{article}
\usepackage{latexsym}
\usepackage{amsmath}
\usepackage{amsfonts}
\usepackage{amssymb}
\renewcommand{\a }{\alpha }
\renewcommand{\b }{\beta }
\renewcommand{\d}{\delta }

\newcommand{\D }{\Delta }

\newcommand{\e }{\varepsilon }

\renewcommand{\i }{\iota}
\renewcommand{\l }{\lambda }

\newcommand{\n }{\nabla }
\newcommand{\var }{\varphi }

\newcommand{\Sig }{\Sigma}

\newcommand{\ov}{\overline}
\newcommand{\wtilde }{\widetilde}

\newcommand{\be}{\begin{equation}}
\newcommand{\ee}{\end{equation}}
\newenvironment{pf}{\noindent{\bf Proof.}\enspace}{%\rule{2mm}{2mm}
\hfill$\Box$\medskip}
\newenvironment{pfn}[1]{\noindent{\bf Proof of {#1}\enspace}}{%\rule{2mm}{2mm}
\hfill$\Box$\medskip}
\newcommand{\R}{\mathbb{R}}
\newtheorem{thm}{Theorem}[section]
\newtheorem{pro}[thm]{Proposition}
\newtheorem{lem}[thm]{Lemma}
\newtheorem{rem}[thm]{Remark}
\newtheorem{cor}[thm]{Corollary}

\numberwithin{equation}{section}

\author{Mohamed BEN AYED$^a$ \& Khalil EL MEHDI$^{b,c}$\thanks{Corresponding author. Present adress : see adress ``c'' mentioned above.}\\ 
{\footnotesize
a : D{\'e}partement de Math{\'e}matiques, Facult{\'e} des Sciences de Sfax, Route
Soukra,}\\
{\footnotesize
 Sfax, Tunisia. E-mail : \texttt{Mohamed.Benayed@fss.rnu.tn} }\\
{\footnotesize
b :  Facult\'e des Sciences et Techniques, Universit\'e de Nouakchott, BP 5026,}\\{\footnotesize
  Nouakchott, Mauritania. E-mail : \texttt{khalil@univ-nkc.mr}}\\
{\footnotesize
 c : The Abdus Salam ICTP, Mathematics Section,  Strada Costiera 11,}\\
{\footnotesize
 34014 Trieste, Italy. E-mail : \texttt{elmehdik@ictp.trieste.it}}
}

\title { \Large \textbf{Existence of Conformal Metrics on Spheres with\\ 
     Prescribed Paneitz Curvature}}

\begin{document}

\date{ }

\maketitle

{\footnotesize
\noindent
{\bf Abstract.}
In this paper we study the problem of prescribing a fourth order
conformal invariant (the Paneitz curvature) on the $n$-sphere, with $n\geq 5$. Using tools from the theory of  critical points at infinity, we provide some topological conditions on the level sets of a given positive function under which we prove the existence of a metric, conformally equivalent to the standard metric, with prescribed Paneitz curvature.

\medskip\noindent\footnotesize{{\bf Mathematics Subject Classification (2000): }\quad
  35J60, 53C21, 58J05, }\\
\noindent
{\bf Key words :} Fourth order conformal invariant,
Paneitz curvature,  Critical
points at infinity.
}

\section{Introduction and the Main Results}
\label{intro}
Given $(M^4,g)$, a smooth $4$-dimensional Riemannian manifold, let $S_g$ be the scalar curvature of $g$ and $Ric_g$ be the Ricci curvature of $g$. In $1983$, Paneitz \cite{P} discovered the following fourth order operator

$$
P_g^4 \var = \D _g^2 \var - \mbox{ div}_g(\frac{2}{3}S_g - 2 \mbox{Ric}_g )d\var.
$$
This operator is conformally invariant in the sense that if $\tilde{g} = e^{2u}g$ is a metric conformally equivalent to $g$, then

$$
P_{\tilde{g}}^4 \var = e^{-4u}P_g^4 (\var ) \qquad \qquad \mbox{ for all } \var\in C^\infty (M),
$$
and it can be seen as a natural extension of the conformal Laplacian on $2$-manifolds. A generalization of $P_g^4$ to higher dimensions has been discovered by Branson \cite{Br2}. Let $(M,g)$ be a smooth compact Riemannian $n$-manifold, $n\geq 5$. The Paneitz operator  $P^n_g$ of \cite{Br2} is  defined by
$$
P^n_g u= \D ^2_gu - div_g(a_n S_g g + b_nRic_g ) du + \frac{n-4}{2} Q^n_g u,
$$
where
$$
a_n=\frac{(n-2)^2 +4}{2(n-1)(n-2)}, \qquad b_n =\frac{-4}{n-2}
$$
$$
Q^n_g = - \frac{1}{2(n-1)} \D _g S_g + \frac{n^3 -4n^2 + 16n-16}{8(n-1)^2(n-2)^2}S^2_g - \frac{2}{(n-2)^2} |Ric_g|^2.
$$

If $\tilde{g}=u^{4/(n-4)}g$ is a metric conformal to $g$, then for all $\var \in C^\infty (M)$ one has 

$$
P^n_g(u\var ) = u^{(n+4)/(n-4)} P^n_{\tilde{g}}(\var )
$$
and

\begin{eqnarray}\label{e:1}
P^n_g(u) =\frac{n-4}{2} Q^n_{\tilde{g}} u^{(n+4)/(n-4)}.
\end{eqnarray}
For more details about the properties of the Paneitz operator, see for example \cite{Br2}, \cite{BCY}, \cite{CGY1}, \cite{CQY1}, \cite{CQY2}, \cite{CY1}, \cite{DHL}, \cite{DMO}, \cite{F}.

In view of relation \eqref{e:1}, it is natural to consider the problem of prescribing the conformal invariant $Q^n$ called the Paneitz curvature, that is :  given a function $ f: M \to \R$, does there exist a metric $\tilde{g}$ conformally equivalent to $g$ such that $Q^n_{\tilde{g}}=f$ ? By equation \eqref{e:1}, the problem can be formulated as follows. We look for solutions of 

\begin{eqnarray}\label{e:2}
P^n_g(u) =\frac{n-4}{2} f u^{(n+4)/(n-4)}, \qquad u > 0\qquad \mbox{on}\qquad M.
\end{eqnarray}

In this paper we consider the case of the standard sphere $(S^n, g)$, $n\geq 5$. Thus,  we are reduced to finding a positive solution $u$ of the problem

\begin{eqnarray}\label{e:3}
 \mathcal{P}u=\D^2u-c_n\D u+d_n u =K u^\frac{n+4}{n-4}, \qquad u > 0\qquad \mbox{ on } S^n,
\end{eqnarray}
where $c_n=\frac{1}{2}(n^2-2n-4)$, $d_n=\frac{n-4}{16}n(n^2-4)$ and where $K$ is a given function defined on $S^n$.

Problem \eqref{e:3} can be viewed as the analogue of the well known scalar curvature problem on $S^n$
\begin{eqnarray}\label{e:cs}
-\D u + \frac{n(n-2)}{4} u = K u^{\frac{n+2}{n-2}},\qquad u>0\qquad \mbox{on}\qquad S^n
\end{eqnarray}
to which many works have been devoted. For details see \cite{A1}, \cite{A}, \cite{AH}, \cite{B2}, \cite{BC2}, \cite{BCCH}, \cite{Bi}, \cite{CGY}, \cite{CL}, \cite{CY2}, \cite{CY3}, \cite{CD}, \cite{ES}, \cite{Ha}, \cite{KW}, \cite{yl}, \cite{Sc1}, \cite{SZ}, \cite{WX} and the references therein.

As for \eqref{e:cs}, there are topological obstructions to solve \eqref{e:3}, based on a Kazdan-Warner type condition, see \cite{DHL} and \cite{WX}. Thus a natural question arises : under which conditions on $K$, does \eqref{e:3} admit a solution? In this paper, we  give sufficient conditions on $K$ such that \eqref{e:3} possesses a solution.

Notice that, problem \eqref{e:cs} has been widely studied in the last decades. On the other hand, to  the authors' knowledge problem \eqref{e:3} has been studied in \cite{BE}, \cite{DHL}, \cite{DMO} and \cite{F} only. In \cite{BE}, \cite{DHL}, \cite{DMO}, the authors treated the lower dimensional case ($n=5,6$). In \cite{F}, Felli proved a perturbative theorem and some existence results under  assumptions of symmetry.

In this paper, we give a contribution in the same direction as in the papers of  Aubin-Bahri \cite{AB} and Bahri \cite{B2} where the problem of prescribing the scalar curvature on closed manifolds was studied using some topological and dynamical tools of the theory of critical points at infinity, see Bahri \cite{B1}. We extend these tools to the framework of such higher order equations. The technique is to use the difference of topology between the level sets of the function $K$ to create a critical point of the Euler functional $J$ associated to \eqref{e:3}. Under our conditions on $K$, the main issue is to prove that the difference in topology between the level sets of $J$ is not completely created by critical points at infinity. This then implies the existence of a critical point of $J$. However, in our situation we have to prove in addition that this critical point is given by a positive function on $S^n$. It is known that in the framework of  higher order equations such a proof is quite difficult in general (see \cite{DMO} for example), and the way we handle this problem here is very simple compared with the literature.

To state our main results, we need to introduce the assumptions that we will  use and some notations.\\
${\bf (A_0)}$\hskip 0.3cm We assume  that $K$ is a positive $C^3$-function on $S^n$ and which has
only nondegenerate critical points $y_0,...,y_l,...,y_s$  with
$$
K(y_0)\geq K(y_1)\geq ...\geq K(y_l)> K(y_{l+1})\geq ...\geq
K(y_s), \,\,\mbox{where}\,\, 0\leq l \leq s.
$$
${\bf (A_1)}$\hskip 0.3cm We assume  that 
$$
-\D K(y_i)>0,\, \mbox{ for } 0\leq i\leq l, \qquad -\D K(y_i)<
0, \,  \mbox{ for } l+1\leq i\leq s.
$$
${\bf (A_1')}$ \hskip 0.3cm We assume that 
$$
 -\D K(y_i)<0,\quad  \mbox{ for }\quad  l+1\leq i\leq s.
$$
In addition, for every $i\in \{1,...,l\}$ such that $-\D K(y_i) \leq 0$, we assume that
$$
n-m+3 \leq \mbox{ index} (K, y_i) \leq n-2,
$$
where index$(K,y_i)$ is the Morse index of $K$ at $y_i$ and $m$ is an integer defined in assumption $(A_2)$.\\
${\bf (A_2)}$ \hskip 0.3cm We assume that there exists a pseudogradient $Z$ of $K$ of
Morse-Smale type (that is, the intersections of the stable and the
unstable manifolds of the critical points of $K$ are transverse) such that the set $ X=\overline{\cup_{0\leq i\leq l}W_s(y_i)}$, where $W_s(y_i)$ is the stable manifold of $y_i$ for $Z$, is not contractible, and we denote by $m$ the
dimension of the first nontrivial reduced homology group of $X$.\\
${\bf (A_3)}$ \hskip 0.3cm We assume that there exists   a
positive constant $\bar{c}$ such that $\bar{c} < K(y_l)$ and such that $X$ is
deformable to a point  in
$K^{\bar{c}}=\{x\in S^n\, \mid \, K(x)\geq \bar{c}\}$.

Now, we are able to state our main results.

\begin{thm}\label{t:11}
There exists  a positive constant $c_0$ independent of $K$
such that if $K$ satisfies $(A_0)$, $(A_1)$,
$(A_2)$, $(A_3)$ and   $K(y_0)/\bar{c}\leq 1+c_0$, then \eqref{e:3} has a solution.
\end{thm}

\begin{cor}\label{c:12}
The solution obtained in Theorem \ref{t:11} has an augmented \\Morse 
index greater than or equal to $m$.
\end{cor}

\begin{thm}\label{t:13}
Assume that $n\geq 6$. Then, there exists  a positive constant $c_0$ independent of $K$
such that if $K$ satisfies $(A_0)$, $(A_1')$,
$(A_2)$, $(A_3)$ and   $K(y_0)/\bar{c}\leq 1+c_0$, then \eqref{e:3} has a solution.
\end{thm}
\begin{rem}
{\bf i).} The assumption $n\geq 6$ in Theorem \ref{t:13} is needed in order to make $(A_1')$ meaningful.\\
{\bf ii).} The assumption $K(y_0)/\bar{c}\leq 1+c_0$ allows basically to perform a single-buble analysis.\\
{\bf iii).} To see how to construct an example of a function $K$ satisfying our assumptions, we refer the interested reader to \cite{AB2}.
\end{rem}
The present paper is organized as follows. In section 2, we set up the variational structure and recall some known facts. In section 3, we prove a Morse lemma at infinity for the Euler functional associated to \eqref{e:3}. In section 4, we provide the proof of Theorem \ref{t:11} and Corollary \ref{c:12}, while section 5 is devoted to the proof of Theorem \ref{t:13}.

\section{Variational Structure and Some Known Facts}
\label{sec:1}

In this section we assume that $K$ is a positive $C^3$-function and we are going to recall the functional setting, the variational problem and its main features.

For $K\equiv 1$, the solutions of \eqref{e:3} form a family
$\wtilde{\d}_{(a,\l)}$ defined by

$$
\wtilde{\d}_{(a,\l)}(x)=\b_n \frac{1}{2^{\frac{n-4}{2}}}
\frac{\l^\frac{n-4}{2}}{\bigl(1+\frac{\l^2-1}{2}(1-\cos
 d(x,a))\bigr)^\frac{n-4}{2}},
$$
 where $a\in S^n$,  $\l >0$ and $\b_n= [(n-4)(n-2)n(n+2)]^{(n-4)/8}$.

 After performing a stereographic projection $\pi$ with the
 point $-a$ as pole, the function $\wtilde{\d}_{(a,\l)}$ is
 transformed into

 $$
\d_{(0,\l)}=\b_n\frac{\l^{\frac{n-4}{2}}}{(1+\l^2\mid
 y\mid^2)^{\frac{n-4}{2}}},
$$
 which is a solution of the problem (see \cite{Li})
$$
 \D^2u= u^\frac{n+4}{n-4} ,\, u>0\,\quad \mbox{ on } \quad \R^n.  
$$

The space $ H_2^2(S^n)$ is equipped with the norm :
$$
\mid\mid u\mid\mid^2= \langle u,u\rangle_\mathcal{P}=\int_{S^n} \mathcal{P}u\cdot u = \int_{S^n} \mid \D u\mid^2 +c_n
\int_{S^n} \mid\n u\mid^2 + d_n\int_{S^n} u^2.
$$
We denote by $\Sig$ the unit sphere of $H_2^2(S^n)$ and we set\\
$\Sig ^+ =\{u\in \Sig \mid u >0\}$.\\
We introduce  the following functional defined on $\Sig$ by
$$
J(u)=\frac{1}{(\int_{S^n}
K|u|^{\frac{2n}{n-4}})^{\frac{n-4}{n}}}=\frac{\mid\mid
u\mid\mid^2}{(\int_{S^n} K|u|^{\frac{2n}{n-4}})^{\frac{n-4}{n}}}.
$$
 The positive critical points of $J$, up to multiplicative constant, are solutions of \eqref{e:3}. The Palais-Smale condition fails for $J$ on $\Sig ^+$. This failure can be described, following the ideas introduced in \cite{BrC}, \cite{L}, \cite{S} as follows :
\begin{pro}\label{p:21}
Assume that $J$ has no critical point in $\Sig ^+$ and
let $(u_k)$ be a sequence in $\Sig ^+$ such that $J(u_k)$ is bounded
and $\n J(u_k)$ tends to 0. Then, there exist an integer $p$ and a
sequence $\e _k$ such that $u_k\in V(p,\e _k)$, where $V(p,\e)$ is
defined by
\begin{align*}
V(p,&\e) =\biggl\{u\in  \Sig \mid \exists a_1, ..., a_p \in S^n, \exists \l _1,..., \l _p > \e^{-1}, \exists \a _1,..., \a _p > 0  \mbox{ with }\\ 
 &\left.\bigg |\bigg | u-\sum_{i=1}^p\a_i\wtilde{\d}_{(a_i, \l_i)}
\bigg |\bigg | < \e ; \, \,   \bigg |\frac{\a _i ^{8/(n-4)}K(a_i)}{\a _j ^{8/(n-4)}K(a_j)}-1\bigg |<\e \,  \forall i, j, \, \,  \e _{ij} < \e \, \forall i\neq j \right\}.
\end{align*}
Here 
$$
\e_{ij}=\left(\frac{\l_i}{\l_j}+\frac{\l_j}{\l_i}+\frac{\l_i\l_j}{2}(1-\cos d(a_i,a_j))\right)^{-\frac{n-4}{2}}.
$$
\end {pro}
We then have the following result which 
defines a parametrization of the set $V(p,\e )$. It follows from corresponding statements in \cite{B2}, \cite{BaC}.
 \begin{pro}\label{p:22} 
For any $p\in N^*$, there exists $\e_p>0$ such that, if
$0<\e<\e_p$ and  $u\in V(p,\e)$, then the following minimization problem
$$
\min \biggl\{\bigg |\bigg | u-\sum_{i=1}^p\a_i\wtilde{\d}_{(a_i, \l_i)}
\bigg |\bigg |,\,  \a_i>0,\, \l_i>0 ,\,  a_i\in S^n\biggr\}
$$
 has a unique solution $(\a,a,\l)=(\a _1,...,\a _p,a_1,...,a_p,\l _1,...,\l _p)  $ (up to permutation).
In particular, we can write $u \in V(p,\e )$ as follows 
$$
u=\sum_{i=1}^p\a_i\wtilde{\d}_{(a_i,\l_i)} +v,
$$
where $v\in H^2_2(S^n)$ such that
\begin{align}\label{V0}
(V_0):\quad
\langle& v,\varphi_i\rangle_\mathcal{P}=0\,\,\mbox{ for }\,\,  i=1,...,p\,\,\mbox{ and every }\\
&\varphi_i = \wtilde{\d}_{(a_i,\l_i)},\,\, \partial\wtilde{\d}_{(a_i,\l_i)}/\partial\l_i,\,\, \partial\wtilde{\d}_{(a_i,\l_i)}/(\partial a_i)_j,\, j=1,..., n,\notag\\
&\mbox{for some system of coordinates }\, (a_i)_1,..., (a_i)_{n}\,\,\mbox{on } S^n\,\,\mbox{near } a_i.\notag    
\end{align}
\end{pro}
Next, we recall from \cite{BE} a useful expansion
 of the functional associated to \eqref{e:3} and its gradient near a critical point at infinity.
 \begin{pro}\label{p:23} \cite{BE}
For $\e$ small enough and $u=\sum_{i=1}^p\a _i\wtilde{\d}_i+v \in V(p,\e)$, the following expansion holds
\begin{align*}
J(u)&= \frac{(\sum_{i=1}^p\a _i ^2) S_n^{4/n} }{(\sum_{i=1}^p\a_i
^{\frac{2n}{n-4}}K(a_i))^{\frac{n-4}{n}}} \biggl[1-\frac{c_2(n-4)}{n}
\sum_{i=1}^p \frac{4\D K(a_i)\a_i ^{\frac{2n}{n-4}}}{\l_i^2\sum_{j=1}^p\a_j
^{\frac{2n}{n-4}}K(a_j)S_n}\\
 & +c_1\sum_{i\ne j}\a
 _i\a_j\e_{ij}\biggl(\frac{1}{\sum_{i=1}^p\a_i ^2
 S_n}-2\frac{\a_i ^{\frac{8}{n-4}}K(a_i)}{\sum_{i=1}^p\a_i
 ^{\frac{2n}{n-4}}K(a_i)S_n}\biggr) -f(v)\\
&+ \frac{1}{\sum_{i=1}^p\a_i ^2
 S_n}Q(v,v) 
 + o\biggl(\sum \e_{kr} + \sum \frac{1}{\l_i^2}\biggr)+ O\biggl( \mid\mid
v\mid\mid^{\min (\frac{2n}{n-4},3)}\biggr) \biggr]
\end{align*}

where $c_1=\b_n^{2n/(n-4)}\int_{\R^n} \frac{dx}{(1+|x|^2)^{(n+4)/2}}$, $c_2 = \frac{1}{2n}\int_{\R^n}|x|^2 \d_{(0,1)}^{2n/(n-4)}$,\\ $S_n=\int_{R^n}\d_{(0,1)}^{2n/(n-4)}$,
$$Q(v,v)=\mid\mid v\mid\mid^2- \frac{n+4}{n-4}\frac{\sum_{i=1}^p\a_i ^2
 }{\sum_{i=1}^p\a_i
^{\frac{2n}{n-4}}K(a_i)}
\int_{S^n}K(\sum_{i=1}^p\a_i\wtilde{\d}_i)^{\frac{8}{n-4}}v^2
$$

and
$$ f(v)=\frac{2}{\sum_{j=1}^p\a_j
 ^{\frac{2n}{n-4}}K(a_j)S_n}
 \int_{S^n}K(\sum_{i=1}^p\a_i\wtilde{\d}_i)^{\frac{n+4}{n-4}}v.
 $$
(Here and in the sequel $\wtilde{\d}_i$ denotes $\wtilde{\d}_{(a_i,\l _i)}$.)
\end{pro}
 Let us introduce the following set 
$$
E=\{v| \, v  \mbox{  satisfies } (V_0) \mbox{ and } ||v|| <\e \},
$$
where $(V_0)$ is defined in \eqref{V0}.
\begin{pro}\label{p:24} \cite{BE}
For any $u=\sum_{i=1}^p\a_i\wtilde{\d}_i \in V(p,\e)$, there
exists a unique $\ov{v}=\ov{v}(a,\a,\l)$ which minimizes $J(u+v)$
with respect to $v\in E$. Moreover, we have the following estimate
\begin{align*}
\mid\mid \ov{v}\mid\mid\leq c\mid\mid f\mid\mid\leq c \biggl[ & 
\sum_{i=1}^p\left(\frac{\mid\n K(a_i)\mid}{\l_i}+\frac{1}{\l_i
  ^2}\right)\\
 & + \sum_{i\ne
  j}\e_{ij}^{\min \left(1, \frac{n+4}{2(n-4)}\right)}(\log\e_{ij}^{-1})^{\min\left(\frac{n-4}{n}, \frac{n+4}{2n}\right)}\biggr].
\end{align*}
\end{pro}
\begin{pro}\label{p:25}\cite{BE}
For $\e>0$ small enough and $u=\a\wtilde{\d}_{(a,\l)}\in V(1,\e)$, the following expansions hold
\begin{align*}
\langle\n J(u),\l\frac{\partial \wtilde{\d}}{\partial \l}\rangle_\mathcal{P} &=
\frac{8(n-4)}{n}c_2\a_i
^{\frac{n+4}{n-4}}J(u)^{\frac{2n-4}{n-4}}\frac{\D K(a_i)}{\l_i ^2}
 +o( \frac{1}{\l_i ^2}),\\
\langle\n J(u),\frac{1}{\l}\frac{\partial \wtilde{\d}}{\partial a}\rangle_\mathcal{P} &= -2c_3
J(u)^{\frac{2n-4}{n-4}}\frac{\n K(a)}{\l}+O(\frac{1}{\l^2}).
\end{align*}
\end{pro}

\section{Morse Lemma at infinity }
\label{sec:2}

In this section, we assume that $K$ is a positive $C^3$-function on $S^n$ having  only nondegenerate critical points $y_0,...,y_l,...,y_s$ satisfying $\D K(y_i) \ne 0$ for every $i \in \{0,...,s\}$, and we consider the case of a single mass. Our goal here is to  perform a Morse lemma at infinity for $J$, which completely gets rid of the $v$-contribution and shows that the functional behaves as $J(\a\wtilde{\d}_{(\wtilde{a},\wtilde{\l})})+ \mid V\mid^2$, where $V$ is a variable independent of $\tilde{a}$ and $\tilde\l$ belonging to a neighborhood of zero in a fixed Hilbert space. Namely, we prove the following result.
\begin{pro}\label{p:31}
For $\e>0$ small enough, there is a diffeomorphism
$$\a\wtilde{\d}_{(a,\l)}+v\longrightarrow (
\a\wtilde{\d}_{(\wtilde{a},\wtilde{\l})},V)$$ such that
\begin{eqnarray}\label{e:p31}
J(\a\wtilde{\d}_{(a,\l)}+v)=J(\a\wtilde{\d}_{(\wtilde{a},\wtilde{\l})})+
\mid V\mid^2.
\end{eqnarray}
Moreover, \eqref{e:p31} can be improved when the concentration point
is near a critical point $y$ of $K$ with $-\D K(y)>0$, leading to
the following normal form :
there is another change of variable
$(\wtilde{a},\wtilde{\l})\longrightarrow
(\bar{a},\bar{\l})$ such that
$$
J(\a\wtilde{\d}_{(\wtilde{a},\wtilde{\l})})=\Psi
(\bar{a},\bar{\l}):=
\frac{S_n^{4/n}}{K(\bar{a})^{(n-4)/n}}\left(1-\frac{4(n-4)}{nS_n}c_2 \frac{(1-\eta)}
{\bar\l^2}\frac{\D K(y)}{K(y)}\right),
$$
where $c_2$ is defined in Proposition \ref{p:23} and where $\eta$ is a small positive constant.
\end{pro}
The proof of Proposition \ref{p:31} can be easily deduced from the following lemma, arguing as in \cite{B2} and \cite{BCCH}.
\begin{lem}\label{l:32}
There exists a pseudogradient $W$ such that the
following holds: \noindent There is a constant $c>0$ independent
of $u=\a\wtilde{\d}_{(a,\l)}
\in V(1,\e)$ such that\\
1) $\langle -\n J(u),W\rangle_\mathcal{P}\geq c(\frac{\mid\n K(a)\mid}{\l}+\frac{1}{\l^2})$\\
2) $\langle -\n J(u+\ov{v}),W+\frac{\partial \ov{v}}{\partial
(\a,a,\l)}(W)\rangle_\mathcal{P} \geq c(\frac{\mid\n K(a)\mid}{\l}+\frac{1}{\l^2})$\\
3) $W$ is bounded\\
4) the only region where $\l$ increases along the flow lines of $W$ is
where $a$ is near a critical point $y$ of $K$ with $-\D K(y)>0$.
\end{lem}

\begin{pf}
Let $\mu > 0$ be such that for any critical point $y$ of $K$, if $d(x,y)\leq 2\mu$ then $\mid\D K(x)\mid > c >0$. Three cases may occur.\\
\noindent {\bf Case 1} $d(a,y)>\mu$ for any critical point $y$. In this case we have $\mid\n K(a)\mid > c >0$. Set
$$Z_1=\frac{1}{\l}\frac{\partial \wtilde{\d}}{\partial
a}\frac{\n K(a)}{\mid \n K(a)\mid}.$$
 From Proposition \ref{p:25}, we have
$$\langle -\n J(u),Z_1\rangle_\mathcal{P}\geq c\frac{\mid\n K(a)\mid}{\l}+O(\frac{1}{\l^2})\geq
c(\frac{\mid\n K(a)\mid}{\l}+\frac{1}{\l^2}).$$
 \noindent
  {\bf Case 2}
  $d(a,y)\leq 2\mu$, where $y$ is a
critical point of $K$ with $-\D K(y)<0$. Set
$$
Z_2=-\l\frac{\partial \wtilde{\d}}{\partial \l}+m_1\var(\l\mid\n K(a)\mid)Z_1,
$$
where $m_1$ is a small constant and $\var$ is a $C^\infty$ function which satisfies $\var(t)=1$ if $t\geq 2$ and $\var(t)=0$ if $t\leq 1$.\\
 Using Proposition \ref{p:25}, we derive 
$$
\langle -\n J(u),Z_2\rangle_\mathcal{P}\geq \frac{c}{\l^2}
+m_1 c (\frac{\mid\n K(a)\mid}{\l}+O(\frac{1}{\l^2})) \geq c(\frac{\mid\n
K(a)\mid}{\l}+\frac{1}{\l^2}).
$$
 \noindent {\bf Case 3} $d(a,y)\leq
2\mu$, where $y$ is a critical point of $K$ with $-\D K(y)>0$. Set
$$
Z_3=\l\frac{\partial \wtilde{\d}}{\partial \l}+m_1
\var(\l\mid\n K(a)\mid)Z_1.
$$ 
We obtain the same equality as in case 2.\\
Hence $W$ will be built as a convex combination of $Z_1$, $Z_2$
and $Z_3$. Thus the proof of claim 1) is completed. Claims 3) and 4) follow easily  from the definition of $W$.
Estimate 2) can be obtained,  arguing as in \cite{B2} and \cite{BCCH}, using Claim 1). 
\end{pf}\\
 Next, we derive from the above results the characterization of the critical points at infinity in $V(1,\eta )$. We recall that  critical points at infinity are the orbits of the gradient flow that remain in $V(p,\e (s))$, where $\e (s)$ is some function such that $\e (s)$ tends to zero when the flow parameter $s$ tends to $+\infty$ (see \cite{B1}).
\begin{pro}\label{p:33}
The only critical points at infinity of $J$ in $V(1,\e )$ correspond to $\tilde\d_{(y,\infty)}$, where $y$ is a critical point of $K$ with $-\D K(y) > 0$.
\end{pro}
\begin{pf}
 From Lemma \ref{l:32}, we know that the only region where $\l$ increases along the pseudogradient $W$, defined in Lemma \ref{l:32}, is the region where $a$ is near a critical point $y$ of $K$ with $-\D K(y) > 0$.
Proposition \ref{p:31} yields a split of variables $a$ and $\l$, thus it is easy to see that if $a=y$, only $\l$ can move. To decrease the functional $J$, we have to increase $\l$, thus we obtain a critical point at infinity only in this case and our result follows.
\end{pf}

\section{Proof of Theorem \ref{t:11} and Corollary \ref{c:12}}
\label{sec:3}
 \begin{pfn}{Theorem \ref{t:11}}
 Arguing by contradiction, we suppose that $J$ has no
critical points in $V_\eta (\Sig^+)$, where
\begin{eqnarray}\label{41}
V_\eta (\Sig ^+)=\{u\in \Sig \mid J(u)^{\frac{2n-4}{n-4}}e^{2J(u)} |u^-|_{L^{2n/(n-4)}}^{8/(n-4)} < \eta\},
\end{eqnarray}
 $\eta$ is a small positive constant and $u^-$ denotes the negative part of $u$, that is, $u^-=\max (0,-u)$.\\
According to Lemma 5.1 in \cite{BE}, we know that $V_\eta (\Sig ^+)$ is invariant under the flow generated by $-\n J$.
 It follows from Proposition \ref{p:33}, that under
the assumptions of Theorem \ref{t:11}, the critical points at
infinity of $J$ under the level $c_1=
(S_n)^{\frac{4}{n}}(K(y_l))^{\frac{4-n}{n}} + \e $ , for $\e$
small enough, are in one to one correspondence with the  critical
points of $K$, i.e. $y_0$, $y_1$, ..., $y_l$. The unstable manifold of such  critical points at infinity, $W_u(y_0)_\infty$,
..., $W_u(y_l)_\infty$, can be described, using Proposition
\ref{p:31}, as the product of $W_s(y_0)$, ..., $W_s(y_l)$ (for a
pseudogradient of $K$ ) by $[A, +\infty [$,
 domain of the variable $\l$, for some positive
number $A$ large enough.\\
Since $J$ has no critical points in $V_\eta (\Sig^+)$, it follows that \\$ J_{c_1}=\{u\in
V_\eta (\Sig ^+) \mid J(u) \leq c_1 \}$ retracts by deformation onto $X_\infty =
\cup _{0\leq j\leq l}W_u(y_j)_\infty$ (see Sections 7 and 8 of
\cite{BR}) which can be parametrized  by $X
\times
[A, +\infty[$.\\
On the other hand, $X_\infty$ is contractible in
$J_{c_2+\e}$, where $c_2=(S_n)^{\frac{4}{n}}\bar{c}^{\frac{4-n}{n}}$.
Indeed, from $(A_3)$, it follows that there exists a continuous contraction $
h :[0,1] \times X \to K^{\bar{c}}$ such that for any $a\in
X \quad h(0,a)=a$ and $h(1,a)=a_0$, a point of $X$. Such a
contraction  gives rise to the following contraction $\tilde{h} :
X_\infty \to V_\eta (\Sig ^+)$ defined by
$$
[0,1] \times X \times \left[A,\right.+\infty\left[ \right.\ni
(t,a,\l  ) \longmapsto \tilde\d _{(h(t,a),\l )} + \bar{v}
\in V_\eta (\Sigma ^+).
$$
For $t=0$, $\tilde\d _{(h(0,a),\l )}+\bar{v} = \tilde\d _{(a,
\l) } +\bar{v} \in X_\infty$. Also, $\tilde{h}$ is continuous and
$\tilde{h}(1,a,\l )= \tilde\d _{(a_0,\l)} +\bar{v}$, hence
our
claim follows.

Now, using Proposition \ref{p:23}, we deduce that
$$
J(\tilde\d _{h(t,a), \l } + \bar{v}) \sim
({S_n})^{\frac{4}{n}}(K(h(t,a)))^{\frac{4-n}{n}}\left(1+O(A ^{-2})\right),
$$
where $K(h(t,a_1)) \geq \bar{c} $ by construction.
Therefore such a contraction is performed below the level $c_2 +\e$ (for $A$
large enough), so $X_\infty$ is contractible in $J_{c_2+\e }$.

In addition, choosing $c_0$ small enough, we see that there is no critical point at infinity for $J$ between the levels $c_2+\e$ and $c_1$, thus  $J_{c_2+\e }$ retracts
by deformation onto $J_{c_1}$, which in turn retracts by deformation onto
$X_\infty$. Therefore $X_\infty$ is contractible leading also to the
contractibility of $X$, which is in contradiction with our
assumption. Hence $J$ has a critical point in $V_\eta (\Sig ^+ )$.

 Now, it remains to prove that such a critical point is a positive function.
Let us define the function $w^-$ by the solution of the following
problem
$$\mathcal{P} w^-=-K(x) (u^-)^{\frac{n+4}{n-4}}\quad \mbox{on }\, S^n.$$
Since $K(x) (u^-)^{\frac{n+4}{n-4}}\in L^{\frac{2n}{n+4}}$, we see $w^-\in H_2^2$. Furthermore, we have $w^-\leq 0$. Thus we derive 
 \begin{align*}
  \int_{S^n} \mathcal{P} w^-\cdot w^- & =\mid\mid w^-\mid\mid^2=\int_{S^n}
 -K(x) (u^-)^{\frac{n+4}{n-4}}w^-\\
&\leq C\mid
 w^-\mid_{L^{2n/(n-4)}}\mid u^-\mid_{L^{2n/(n-4)}}^{(n+4)/(n-4)}\\
  & \leq C\mid\mid w^-\mid\mid \, \mid u^-\mid_{L^{2n/(n-4)}}^{(n+4)/(n-4)}.
  \end{align*}
 Consequently, either $||w^-||=0$ and therefore $u^-=0$, or $||w^-||\ne 0$ and we derive 
 \begin{eqnarray}\label{w-1}
 \mid\mid w^-\mid\mid\leq C\mid
 u^-\mid_{L^{2n/(n-4)}}^{(n+4)/(n-4)}.
\end{eqnarray}
Furthermore, on  one hand we have
\begin{align}\label{w-2}
 \int_{S^n}u\cdot \mathcal{P} w^-=\int_{S^n}-u&
 K(u^-)^{\frac{n+4}{n-4}}=\int_{S^n}K(u^-)^{\frac{2n}{n-4}}\notag\\
&\geq
 c_K\int_{S^n}(u^-)^{\frac{2n}{n-4}}\geq c_K \mid
 u^-\mid_{L^{\frac{2n}{n-4}}}^{\frac{2n}{n-4}}
 \end{align}
 (since $K$ is bounded from below by a positive constant $c_K$), and on
  the other hand, we have
 \begin{align}\label{w-3}
 \int u\cdot\mathcal{P}w^- & =\int w^-\cdot\mathcal{P} u=\int w^- K\mid u\mid^{\frac{8}{n-4}}u\leq 
 \int_{u\leq 0}-w^-K(u^-)^{\frac{n+4}{n-4}}\notag \\
  & \leq \int_{S^n}-w^-K(u^-)^{\frac{n+4}{n-4}}=\int_{S^n}
  w^-\cdot\mathcal{P}w^-=\mid\mid w^-\mid\mid^2.
  \end{align}
  Using \eqref{w-1}, \eqref{w-2} and \eqref{w-3}, we obtain
 $$c_K\mid u^-\mid_{L^{2n/(n-4)}}^{2n/(n-4)}\leq \mid\mid
 w^-\mid\mid^2 \leq C \mid
 u^-\mid_{L^{2n/(n-4)}}^{2(n+4)/(n-4)}.$$
 Observe that $2n/(n-4)<2(n+4)/(n-4)$. Thus, either $u^-=0$ or \\$\mid u^-\mid_{L^{2n/(n-4)}}\geq C$ and
this case cannot occur since by the definition of the
neighborhood of $\Sig ^+$ this norm is small. This completes the proof of our result.
 \end{pfn}\\
\begin{pfn}{\bf Corollary \ref{c:12}}
Arguing by contradiction, we may assume that the Morse index of
the
solution provided by Theorem \ref{t:11} is not exceeding  $m-1$.\\
Perturbing $J$, if necessary, we may assume that all the critical
points of $J$ are nondegenerate and have their Morse index less than $m-1$. Such critical points do not change the homology group in
dimension $m$ of the level sets of $J$.\\
Since $X_\infty$ defines a homology class in dimension $m$
which is nontrivial in $J_{c_1}$, but trivial in $J_{c_2+\e}$, our
result follows.
\end{pfn}
\section{Proof of Theorem \ref{t:13}}
\label{sec:4}

 Arguing by contradiction, we suppose that $J$ has no critical points in the set $V_\eta (\Sig^+)$ defined by \eqref{41}. Let $\{z_1,...,z_r\}\subset \{y_1,...,y_l\}$  be the critical points of $K$ with 
$$
-\D K(z_j)\leq 0 \quad (1\leq j\leq r).
$$

The idea of the Proof of Theorem \ref{t:13} is to perturb the function $K$ in the $C^1$ sense in some  neigborhoods of $z_1,...,z_r$ such that the new function
$\tilde{K}$ has the same critical points with the same  Morse indices but satisfying that $-\D \tilde{K}(z_j) > 0$ for $1\leq j\leq r$.

The new $\tilde{X}$ corresponding to $\tilde{K}$, defined in 
assumption $(A_2)$, is also not contractible and its homology group in dimension $m$ is nontrivial.

 Under the level $2^{4/n}S_n^{4/n}(K(y_0))^{(4-n)/n}$, the associated functional $\tilde{J}$ is close to the functional $J$ in the $C^1$ sense. Under the level $c_2 +\e$, where $c_2$ is defined in the proof of Theorem \ref{t:11}, the functional  $\tilde{J}$ may have other
critical points, however a careful choice of $\tilde{K}$ ensures that all these critical
points have morse indices less than $m - 2$ (see Proposition \ref{p:51} below), and so they do not change the homology in dimensin $m$, therefore the arguments
used in the Proof of Theorem \ref{t:11} lead to a contradiction. It follows that Theorem \ref{t:13} will be a corollary of the following Proposition:
\begin{pro}\label{p:51}
We can choose $\tilde{K}$ close to $K$ in the $C^1$ sense such that $\tilde{K}$ has the same critical points with the same  Morse indices and such that:
\begin{align*}
i)\quad & -\D \tilde{K}(z_j) > 0 \quad\mbox{ for}\quad  1\leq j\leq r,\\
ii)\quad &  -\D \tilde{K}(y_i) > 0 \quad\mbox{ for}\quad  i\in\{0,...,l\} \diagdown \{1,..., r\},\\
iii)\quad & -\D \tilde{K}(y_i) < 0 \quad\mbox{ for}\quad  l+1\leq i\leq s,\\
iv)\quad &\mbox{ if } \tilde{J} \mbox{ has critical points under the level } c_2 +\e, \mbox{ then their Morse}\\
\quad &\mbox{  indices are less than } m-2, \mbox{ where } m  \mbox{ is defined in assumption } (A_2),\\
v)\quad & \mbox{ the new } \tilde{X} \mbox{ corresponding to } \tilde{K}, \mbox{ defined in assumption }(A_2), \mbox{ is also}\\
\quad & \mbox{  not contractible and its homology group in dimension } m \mbox{ is}\\
\quad & \quad{ nontrivial}.
\end{align*} 
\end{pro}

To prove Proposition \ref{p:51}, we need the following lemmas.
\begin{lem}\label{l:52}
Assume that $n\geq 5$ and let $P=P(z,\l)$ be the $\mathcal{P}$-orthogonal projection onto the linear space generated by  $\tilde{\d}_{(z.\l )}$, $ \partial\tilde{\d}_{(z,\l )}/\partial\l$  and \\$\partial\tilde{\d}_{(z,\l )}/\partial z_i$, $i=1,..., n$, for some system of coordinates $z_1,...,z_{n}$ on $S^n$ near $z$. Then, we have the following estimates
\begin{align*}
&(i)\quad ||J'(\tilde{\d}_{(z.\l )})|| = O\left(\frac{1}{\l}\right); \quad (ii)\quad || \frac{\partial P}{\partial z}||= O(\l);\\
& (iii)\quad || \frac{\partial ^2 P}{\partial ^2 z}||= O(\l ^2).
\end{align*}
\end{lem}
\begin{pf}
The proof of claim $(i)$ is easy, so we will omit it. Now we are going to prove claim $(ii)$. Let
 $$
\var \in \{ \l^{-1}\partial\tilde{\d}_{(z,\l )}/\partial z_i,\,i=1,...,n, \tilde{\d}_{(z.\l )}, \l \partial\tilde{\d}_{(z,\l )}/\partial\l\}.
$$
 We then have $P\var =\var$, therefore
$$
\frac{\partial P}{\partial z_j}(\var )= \frac{\partial \var}{\partial z_j} - P\frac{\partial \var}{\partial z_j},\quad\mbox{for }\,j=1,...,n,
$$
thus $||\frac{\partial P}{\partial z}(\var )||= O(\l||\var||)$.
Now, for $v\in E$, we have $Pv=0$, thus 
$$
\frac{\partial P}{\partial z_j}v = - P \frac{\partial v}{\partial z_j}= \sum_{i=1}^{n+2}a_{ij}\var _i,
$$
where $\var _i=  \l^{-1}(\partial\tilde{\d}_{(z,\l )})/(\partial z_i)$, $i=1,..., n$,   $\var _{n+1} =\tilde{\d}_{(z.\l )}$,  and $\var _{n+2}= \l (\partial\tilde{\d}_{(z,\l )})/(\partial\l)$.\\
But, we have 
$$
a_{ij}||\var _i||^2 = - \langle\frac{\partial v}{\partial z_j}, \var _i\rangle=  \langle v, \frac{\partial\var_i}{\partial z_j}\rangle = O(\l ||v||).
$$
Thus claim $(ii)$ follows.

 In the same way, we obtain claim $(iii)$ and hence the proof of our lemma is completed.
\end{pf}
\begin{lem}\label{l:53}
Let $z_0$ be a point of $S^n$ close to a critical point of $K$ and let $\bar{v}=\bar{v}(z_0,\a , \l )\in E$ defined in Proposition \ref{p:24}. Then, we have the following estimates
$$
(i) \qquad ||\bar{v}||= o(\frac{1}{\l}), \qquad \quad (ii)\qquad ||\frac{\partial \bar{v}}{\partial z}|| = o(1).
$$
\end{lem}
\begin{pf}
We notice that Claim $(i)$ follows from Proposition \ref{p:24}. Then we only need to show that Claim $(ii)$ is true. We know that $\bar{v}$ satisfies
$$
A\bar{v}= f + O\left(||\bar{v}||^{\frac{n+4}{n-4}}\right)\quad \mbox{and } \quad
\frac{\partial A}{\partial z} \bar{v} +A \frac{\partial \bar{v}}{\partial z} = \frac{\partial f}{\partial z} + O\biggl(||\bar{v}||^{\frac{8}{n-4}}\frac{\partial \bar{v}}{\partial z}\biggr),
$$
where $A$ is the operator associated to the quadratic form $Q$ defined on $E$ ($Q$ and $f$ are defined in Proposition \ref{p:23}).
Then we have
$$
A\biggl(\frac{\partial\bar{v}}{\partial z} -P\frac{\partial \bar{v}}{\partial z}\biggr) = \frac{\partial f}{\partial z} - \frac{\partial A}{\partial z} \bar{v}- AP\frac{\partial \bar{v}}{\partial z} +  O\biggl(||\bar{v}||^{8/(n-4)}\frac{\partial \bar{v}}{\partial z}\biggr).
$$
Since $Q$ is a positive quadratic form on $E$ (see \cite{BE}), we then derive 
$$
||\frac{\partial\bar{v}}{\partial z} -P\frac{\partial \bar{v}}{\partial z}|| \leq C\biggl(||\frac{\partial f}{\partial z}|| + ||\frac{\partial A}{\partial z}|||| \bar{v}|| +||P\frac{\partial \bar{v}}{\partial z}||+ || \bar{v}||^{\frac{8}{n-4}}||\frac{\partial \bar{v}}{\partial z}||\biggr).
$$

Now we estimate each term on the right-hand side of the above inequality. First, it is easy to see $||\frac{\partial A}{\partial z}||=O(\l )$. Therefore, by Claim (i), we obtain $||\frac{\partial A}{\partial z}|||| \bar{v}|| =o(1)$. Secondly, we have
\begin{align}
<\frac{\partial f}{\partial z},v> & =c\int K \d_{(z_0,\l)}^{\frac{8}{n-4}}\frac{\partial \d}{\partial z}v\notag\\
&=c\n K(z_0)\int d(z_0,x)\d ^{\frac{8}{n-4}}\frac{\partial \d }{\partial z}v+O\biggl(\int d^2(x,z_0)\d^{\frac{n+4}{n-4}}\l |v|\biggr)\notag \\
 & \leq c ||v||(|\n K(z_0)|+\frac{1}{\l}).
\end{align}
Sinze $z_0$ is close to a critical point of $K$, we derive $||\frac{\partial f}{\partial z}||=o(1)$.
For the  term $||P(\frac{\partial \ov{v}}{\partial z})||$ we have, since $\ov{v}\in E$,
\begin{align*}
<\frac{\partial \ov{v}}{\partial z},\d_{(z_0,\l)}>= & -  <\ov{v},\frac{\partial \d_{(z_0,\l)}}{\partial z}>=0,\\
<\frac{\partial \ov{v}}{\partial z},\l\frac{ \partial \d_{(z_0,\l)}}{{\partial \l}}>= & -  <\ov{v},\l\frac{\partial^2 \d_{(z_0,\l)}}{\partial\l\partial z}>=O(\l||\ov{v}||)=o(1).
\end{align*}
In the same way, we obtain
$$
<\frac{\partial \ov{v}}{\partial z},\frac{1}{\l}\frac{\partial \d_{(z_0,\l)}}{{\partial z}}>=o(1).
$$
Therefore $||P(\frac{\partial \ov{v}}{\partial z})||=o(1)$. Now, using the  inequality 
$$
||\frac{\partial \ov{v}}{\partial z}||\leq ||\frac{\partial \ov{v}}{\partial z} -P(\frac{\partial \ov{v}}{\partial z})|| + ||P(\frac{\partial \ov{v}}{\partial z})||,
$$
we easily derive our claim and the lemma follows.
\end{pf}\\
We are now able to prove Proposition \ref{p:51}.\\
\begin{pfn}{Proposition \ref{p:51}}
We suppose that $J$ has no critical points in $V_\eta (\Sig^+)$ and we perturb  the function $K$ only in some neighborhoods of $z_1,...,z_r$, therefore Claims $ii)$ and $iii)$ follow from assumption $(A_1')$. We notice that under the level $c_2+\e$ and outside $V(1,\e_0)$, we have $|\n J|> c > 0$. If $\wtilde K$ is close to $K$ in the $C^1$ sense, then  $\wtilde J$ is close to $J$ in the $C^1$ sense, and therefore $|\n \wtilde J|> c/2$ in this region. Thus, a critical point $u_0$ of $\wtilde J$ under the level $c_2+\e$ has to be in $V(1,\e_0)$. Thus, we can write $u_0=\wtilde\d_{(z_0,\l)}+\ov{v}$. Now, using Proposition \ref{p:25}, we derive (with a constant $c>0$) 
$$
0=<\n \wtilde J (u_0),\frac{1}{\l}\frac{\partial \d}{\partial z}>=-c\frac{\n \wtilde K(z_0)}{\l}+o(\frac{1}{\l}),
$$
thus $z_0$ has to be close to $y_i$ where $i\in \{0,...,s\}$. We also have, again by Proposition \ref{p:25}, 
\begin{eqnarray}\label{e:51}
0=<\n \wtilde J (u_0),\l\frac{\partial \d}{\partial \l}>=c\frac{\D \wtilde K(z_0)}{\l^2}+o(\frac{1}{\l^2}).
\end{eqnarray} 
In a neighborhood of $y_i$ with $i\in\{k \mid -\D K(y_k) >0\}\cup\{l+1,...,s\}$, we have $\wtilde K\equiv K$ and therefore 
$|\D\wtilde K| >c >0$ in this neighborhood. Thus \eqref{e:51} implies that $z_0$ has to be near $z_i$ with $1\leq i\leq r$ (recall that $z_i$'s are the critical points among $y_1,...,y_l$ with a nonnegative value of  $\D K$).

  In the sequel, we assume that $\wtilde\d=\wtilde\d_{(z,\l)}$ satisfies $||\wtilde\d||=1$, and thus $\mathcal{P}\wtilde\d=S_n^{\frac{4}{n-4}}\wtilde\d^{\frac{n+4}{n-4}}$. We also assume that $|D^2 \wtilde K| \leq c(1+|D^2 K|)$, where $c$ is a fixed positive constant. 

Let $u_0=\wtilde\d_{(z_0,\l)}+\ov{v}$ be a critical point of $\wtilde J$.
In order to compute the Morse index of $\wtilde J$ at $u_0$, we need to compute  $\frac{\partial^2}{\partial z^2}\wtilde J(\wtilde\d_{(z,\l)}+\ov{v})_{|z=z_0}$. We observe that 
 $$\frac{\partial}{\partial z}\wtilde J(\wtilde\d_{(z,\l)}+\ov{v}) = \wtilde J'(\wtilde\d_{(z,\l)}+\ov{v})\frac{\partial}{\partial z}(\wtilde\d_{(z,\l)}+\ov{v}) =
 \wtilde J'(\wtilde\d_{(z,\l)}+\ov{v})P(\frac{\partial}{\partial z}(\wtilde\d_{(z,\l)}+\ov{v}))
$$
and 
\begin{align}\label{e:52} 
\frac{\partial^2}{\partial z^2}\wtilde J(\wtilde\d_{(z,\l)}+\ov{v}) = & \wtilde J''(\wtilde\d_{(z,\l)}+\ov{v})\frac{\partial}{\partial z}(\wtilde\d_{(z,\l)}+\ov{v})P(\frac{\partial}{\partial z}(\wtilde\d_{(z,\l)}+\ov{v}))\\
 & +
 \wtilde J'(\wtilde\d_{(z,\l)}+\ov{v})\frac{\partial}{\partial z}\biggl(P(\frac{\partial}{\partial z}(\wtilde\d_{(z,\l)}+\ov{v}))\biggr).\notag
\end{align}
For $z=z_0$, we have $\wtilde J'(\wtilde\d_{(z,\l)}+\ov{v})=0$. We will estimate each term on the right-hand side of \eqref{e:52}.

 First, we have by Lemma \ref{l:53}
$$
\wtilde J''(\wtilde\d_{(z,\l)}+\ov{v})\frac{\partial \ov{v}}{\partial z}P(\frac{\partial \ov{v}}{\partial z})=o(1).
$$
Secondly, we compute 
\begin{align*}
T&=\wtilde J''(\wtilde\d_{(z,\l)}+\ov{v})\frac{\partial \wtilde \d}{\partial z}P\frac{\partial \ov{v}}{\partial z}\\
&=c\bigg[<\frac{\partial \wtilde \d}{\partial z},P\frac{\partial \ov{v}}{\partial z}>-\frac{n+4}{n-4}\wtilde J(u_0)^{\frac{n}{n-4}}\int \tilde{K}(\wtilde\d+\ov{v})^{\frac{8}{n-4}}\frac{\partial \wtilde \d}{\partial z}P\frac{\partial \ov{v}}{\partial z}\bigg].
\end{align*}
According to Proposition \ref{p:23}, we have 
$$
\wtilde J(\wtilde \d+\ov{v})=\frac{S_n^{4/n}}{\wtilde K(z)^{\frac{n-4}{n}}}+O(\frac{||\ov{v}||}{\l}+\frac{1}{\l ^2}).
$$
Thus 
\begin{align*}
T & = c\left[<\frac{\partial \wtilde \d}{\partial z},P\frac{\partial \ov{v}}{\partial z}> + o(1)\right.\\
&\left. -\frac{n+4}{n-4}S_n^{\frac{4}{n-4}}\int \frac{\tilde{K}}{\tilde{K}(z)}\biggl(\wtilde\d ^{\frac{8}{n-4}}+O\biggl(\wtilde\d^{\frac{12-n}{n-4}}|\ov{v}|+|\ov{v}|^{\frac{8}{n-4}}\chi_{\wtilde\d\leq |\ov{v}|}\biggr)\biggr)\frac{\partial \wtilde \d}{\partial z}P\frac{\partial \ov{v}}{\partial z}\right]\\
& = c\frac{n+4}{n-4}S_n^{\frac{4}{n-4}}\int \left(1-\frac{\tilde{K}}{\tilde{K}(z)}\right)\wtilde\d ^{\frac{8}{n-4}}\frac{\partial \wtilde \d}{\partial z}P(\frac{\partial \ov{v}}{\partial z})\\
&+ O\biggl(\l||\ov{v}||||\frac{\partial \ov{v}}{\partial z}||+\l||\ov{v}||^{\frac{n+4}{n-4}}||\frac{\partial \ov{v}}{\partial z}||\biggr) + o(1)\\
&= o(1).
\end{align*}
Thus \eqref{e:52} becomes 
\begin{align*} 
\frac{\partial^2}{\partial z^2}\wtilde J(&\wtilde\d_{(z,\l)}  +\ov{v}) = \wtilde J''(\wtilde\d_{(z,\l)}+\ov{v})\frac{\partial\wtilde\d}{\partial z}\frac{\partial\wtilde\d}{\partial z}+\wtilde J''(\wtilde\d_{(z,\l)}+\ov{v})\frac{\partial\ov{v}}{\partial z}\frac{\partial\wtilde\d}{\partial z}+o(1)\\
 & =  \wtilde J''(\wtilde\d_{(z,\l)}+\ov{v})\frac{\partial\wtilde\d}{\partial z}\frac{\partial\wtilde\d}{\partial z}+\wtilde J''(\wtilde\d_{(z,\l)})\frac{\partial\ov{v}}{\partial z}\frac{\partial\wtilde\d}{\partial z} + o(1)\\ 
& = \wtilde J''(\wtilde\d_{(z,\l)})\frac{\partial\wtilde\d}{\partial z}\frac{\partial\wtilde\d}{\partial z}+\wtilde J^{(3)}(\wtilde\d_{(z,\l)})\ov{v}\frac{\partial\wtilde\d}{\partial z}\frac{\partial\wtilde\d}{\partial z}
+ \wtilde J''(\wtilde\d_{(z,\l)})\frac{\partial\ov{v}}{\partial z}\frac{\partial\wtilde\d}{\partial z}+o(1),
 \end{align*}
and 
\begin{align*}
\frac{\partial^2}{\partial z^2}\wtilde J(\wtilde\d_{(z,\l)}+\ov{v})& = \wtilde J''(\wtilde\d_{(z,\l)})\frac{\partial\wtilde\d}{\partial z}\frac{\partial\wtilde\d}{\partial z}+\frac{\partial}{\partial z}( \wtilde J''(\wtilde\d_{(z,\l)})\ov{v}\frac{\partial\wtilde\d}{\partial z})\\
&- \wtilde J''(\wtilde\d_{(z,\l)})\ov{v}\frac{\partial^2\wtilde\d}{\partial z^2}+o(1).
\end{align*}
Since 
$$
0=\wtilde J'(\wtilde\d+\ov{v})\frac{\partial^2\wtilde\d}{\partial z^2}=\wtilde J'(\wtilde\d)\frac{\partial^2\wtilde\d}{\partial z^2}+\wtilde J''(\wtilde\d)\ov{v}\frac{\partial^2\wtilde\d}{\partial z^2}+o(1),
$$
we arrive at
\begin{align*} 
\frac{\partial^2}{\partial z^2}\wtilde J(\wtilde\d_{(z,\l)} +\ov{v})& = \wtilde J''(\wtilde\d_{(z,\l)})\frac{\partial\wtilde\d}{\partial z}\frac{\partial\wtilde\d}{\partial z}
+\frac{\partial}{\partial z}( \wtilde J''(\wtilde\d_{(z,\l)})\ov{v}\frac{\partial\wtilde\d}{\partial z})\\
&+\wtilde J'(\wtilde\d)\frac{\partial^2\wtilde\d}{\partial z^2} +  o(1)\\
 & =\frac{\partial^2}{\partial z^2}(\wtilde J(\wtilde\d))+\frac{\partial}{\partial z}
(\wtilde J''(\wtilde\d)\ov{v}\frac{\partial\wtilde\d}{\partial z})+o(1). 
\end{align*}

 Now we notice that 
$$
\wtilde J''(\wtilde\d)\ov{v}\frac{\partial\wtilde\d}{\partial z}=-c\frac{n+4}{n-4}S_n^{\frac{4}{n-4}}\int \biggl(\frac{\wtilde K}{\wtilde K(z)}-1\biggr)\wtilde\d ^{\frac{8}{n-4}}\ov{v}\frac{\partial \wtilde\d}{\partial z}.
$$
Thus 
\begin{align*}
\frac{\partial}{\partial z}(\wtilde J''(&\wtilde\d)\ov{v}\frac{\partial\wtilde\d}{\partial z})  =-c\frac{n+4}{n-4}S_n^{\frac{4}{n-4}}\bigg[\int \frac{-\wtilde K}{\wtilde K(z)^2}\n \wtilde K(z)\wtilde\d ^{\frac{8}{n-4}}\ov{v}\frac{\partial \wtilde\d}{\partial z}\\
& +\int\biggl(\frac{\wtilde K}{\wtilde K(z)}-1\biggr)\wtilde\d ^{\frac{8}{n-4}}\ov{v}\frac{\partial^2 \wtilde\d}{\partial z^2}    
  +\frac{8}{n-4}\int \biggl(\frac{\wtilde K}{\wtilde K(z)}-1\biggr)\wtilde\d ^{\frac{12-n}{n-4}}\ov{v}(\frac{\partial \wtilde\d}{\partial z})^2\\
&+\int \biggl(\frac{\wtilde K}{\wtilde K(z)}-1\biggr)\wtilde\d ^{\frac{8}{n-4}}\frac{\partial\ov{v}}{\partial z}\frac{\partial \wtilde\d}{\partial z}\bigg]\\
 & =O\bigl(\l||\ov{v}||+||\frac{\partial\ov{v}}{\partial z}||\bigr)=o(1).
\end{align*}
Collecting all the estimates, we finally obtain 
$$
\frac{\partial^2}{\partial z^2}\wtilde J(\wtilde\d_{(z,\l)}+\ov{v})=
\frac{\partial^2}{\partial z^2}\wtilde J(\wtilde\d_{(z,\l)})+o(1).
$$
Now we recall that 
$$
\wtilde J(\wtilde\d_{(z,\l)})=\frac{1}{\biggl(\int \wtilde K(x+z)\wtilde\d_{(o,\l)}^{\frac{2n}{n-4}}\biggr)^{\frac{n-4}{n}}}.
$$
Thus 
$$
\frac{\partial}{\partial z}\wtilde J(\wtilde\d_{(z,\l)})=-\frac{n-4}{n}
\frac{\int\n \wtilde K(x+z)\wtilde\d_{(o,\l)}^{\frac{2n}{n-4}}}
{\biggl(\int \wtilde K(x+z)\wtilde\d_{(o,\l)}^{\frac{2n}{n-4}}\biggr)^{\frac{2n-4}{n}}},
$$
 and 
\begin{align*}
\frac{\partial^2}{\partial z^2}\wtilde J(\wtilde\d_{(z,\l)})&=\frac{n-4}{n}\frac{2n-4}{n}
\frac{\biggl(\int\n \wtilde K(x+z)\wtilde\d_{(o,\l)}^{\frac{2n}{n-4}}\biggr)^2}
{\biggl(\int \wtilde K(x+z)\wtilde\d_{(o,\l)}^{\frac{2n}{n-4}}\biggr)^{\frac{3n-4}{n}}}\\
&-\frac{n-4}{n}
\frac{\int D^2 \wtilde K(x+z)\wtilde\d_{(o,\l)}^{\frac{2n}{n-4}}}
{\biggl(\int \wtilde K(x+z)\wtilde\d_{(o,\l)}^{\frac{2n}{n-4}}\biggr)^{\frac{2n-4}{n}}}.
\end{align*}
Observe that 
\begin{align*}
 \int \wtilde K(x+z_0)\wtilde\d_{(o,\l)}^{\frac{2n}{n-4}} & =c\wtilde K(z_0)\,\,\,\,+\,\,O(\frac{1}{\l}),\\
  \int \n \wtilde K(x+z_0)\wtilde\d_{(o,\l)}^{\frac{2n}{n-4}} & =c\n \wtilde K(z_0)\,\,+\,O(\frac{1}{\l}),\\
  \int D^2\wtilde K(x+z_0)\wtilde\d_{(o,\l)}^{\frac{2n}{n-4}} & =cD^2\wtilde K(z_0)+O(\frac{1}{\l}).
\end{align*}
Thus, if $z_0$ is close to a critical point, we have  
$$
\frac{\partial^2}{\partial z^2}\wtilde J(\wtilde\d_{(z,\l)})_{/z=z_0}=-cD^2 \wtilde K(z_0)+o(1).
$$
Without loss of generality, we can assume that $z_0$ is close to $z_1$. We can also assume that $D^2 K=D^2 K(z_1)+o(1)$ in $B(z_1,\rho)$ and $D^2 K(z_1)$ is diagonal, where $\rho$ is a small fixed positive constant. Notice that, by assumption $(A_1')$, $D^2K(z_1)$ admits negative eigenvalues in $B(z_1,\rho)$. Using the diagonal form of $D^2K(z_1)$, we can obtain $\wtilde K$, if we decrease the negative eigenvalues of $D^2 K(z_1)$ in $B(z_1,\rho)$ such that $-\D \wtilde K(z_1) > 0$ and that $\wtilde K$ has only $z_1$ as a critical point in $B(z_1,\rho)$. In this construction, we will bring back the negative eigenvalues of $D^2 K(z_1)$ to their initial values on $\partial B(z_1,\rho)$. The Morse index of $K$ at $z_1$ satisfies index$(K,z_1)\geq n-m+3$.

Since $\rho$ is fixed, the Morse index of $K$ at $z_0$ is equal to the number of negative eigenvalues of $-D^2 \wtilde K(z_0)$ which is the same as the number of negative eigenvalues of $-D^2 K(z_1)$. Thus, the contribution of the variable $z$ to the Morse index of $\tilde{J}$ is less than or equal to $m-3$. Taking into account the contribution of the variable $\l$, we derive  Claims $iv)$ and $i)$. 

On the other hand, according to  assumption $(A_1')$ we have
$$
n-m+3 \leq \mbox{index} (K,z_j)= \mbox{index} (\tilde{K}, z_j)\quad\mbox{for}\quad 1\leq j\leq r.
$$
Thus, for any pseudogradient of $\tilde{K}$, the dimension of the stable manifold of $z_j$ is less than $m-3$. Note that our perturbation changes the pseudogradient $Z$ to $\tilde{Z}$, but only in some neighborhoods of $z_1,..., z_r$. Therefore the stable manifolds of $y_i$, for $\i \not\in \{1,...,r\}$, remain unchanged. Since the dimension of $X$ is greater than $m$ and its homology group in dimension $m$ is nontrivial, we derive that the homology group of $\tilde{X}$ in dimension $m$ is also nontrivial. This completes the proof of Proposition \ref{p:51}. 
\end{pfn}

{\bf Acknowledgements.} Part of this work was done while the first author (M. Ben Ayed) was visiting the Abdus Salam International Centre for Theoretical Physics (Trieste, Italy). He would like to thank this institution for its warm hospitality.

\end{document}